\documentclass[11pt]{amsart}
\usepackage{amssymb, graphicx}
\usepackage[latin1]{inputenc}
\usepackage{amsfonts}
\usepackage{latexsym}
\usepackage{amscd}
\usepackage{enumerate}
\newtheorem{prop}{Proposition}[section]
\newtheorem{lem}[prop]{Lemma}
\newtheorem{thm}[prop]{Theorem}

\newtheorem{ejs}[prop]{Examples}

\theoremstyle{definition}

\numberwithin{equation}{section}

\newcommand{\la}{\langle}
\newcommand{\ra}{\rangle}

\newcommand{\ad}{\mathrm{ad\ }}
\newcommand{\QAn}{\mathrm{QAnn}}
\newcommand{\Ann}{\mathrm{Ann}}
\newcommand{\inn}{\mathrm{Inn}}
\newcommand{\An}{\mathrm{Ann}}
\newcommand{\sder}{\mathrm{SDer}}
\newcommand{\Ad}{\mathrm{ad}\,}

\def\wqa{weak algebra of quotients\ }

\newcommand{\der}{\mathcal{D}\mathrm{er}}

\begin{document}

\title{Strongly non-degenerate Lie algebras}

\author{Francesc Perera \& Mercedes Siles Molina}
\address{Departament de Matem\`atiques,
Universitat Aut\`onoma de Barcelona, 08193 Be\-lla\-terra,
Barcelona, Spain} \email{perera@mat.uab.cat}
\address{Departamento de \'Algebra, Geometr\'\i a y Topolog\'\i a,
Universidad de M\'alaga, 29071 M\'alaga, Spain}
\email{mercedes@agt.cie.uma.es}


\thanks{The first author has been partially supported by the DGI MEC-FEDER through Project MTM2005-00934,
and by the Comissionat per Universitats i Recerca de la Generalitat de Catalunya. The second author has been
partially supported by the MEC and Fondos FEDER jointly through project MTM2004-06580-C02-02 and by the Junta
de Andaluc\'{\i}a PAI, projects FQM-336 and FQM-1215.}





\begin{abstract}
Let $A$ be a semiprime $2$ and $3$-torsion free non-commutative
associative algebra. We show that the Lie algebra $\der(A)$ of
(associative) derivations of $A$ is strongly non-degenerate, which
is a strong form of semiprimeness for Lie algebras, under some
additional restrictions on the center of $A$. This result follows
from a description of the quadratic annihilator of a general Lie
algebra inside appropriate Lie overalgebras. Similar results are
obtained for an associative algebra $A$ with involution and the
Lie algebra $\sder(A)$ of involution preserving derivations of
$A$.
\end{abstract}

\maketitle

\section*{Introduction}

This paper is concerned with the structure of the Lie algebra
$\der(A)$ of associative derivations of an associative algebra
$A$, which we will also assume to be $2$ and $3$-torsion free. It
was proved in~\cite[Theorem 4 and Theorem 2]{JJ} that if $A$ is
semiprime (respectively prime) then $\der (A)$ is a semiprime
(respectively, prime) Lie algebra. We prove below  that, if $A$ is
prime, this result can be strengthened to show that in fact $\der
(A)$ is strongly non-degenerate (see below for the precise
definitions).

The key result in the paper is Theorem~\ref{qadann} and has a
technical flavour. Let $L$ be a subalgebra of a Lie algebra $Q$.
The quadratic annihilator of $L$ inside $Q$ is defined as the set
$\{q\in Q \ \vert\ [q, [q, L]]=0\}$. Roughly speaking,
Theorem~\ref{qadann} allows to obtain non-zero elements in the
quadratic annihilator of $L$ in itself from non-zero elements in
the quadratic annihilator of $L$ inside $Q$ whenever $Q$ is a weak
quotient algebra of $L$, i.e., $[L, q]\neq 0$ for every non-zero
$q\in Q$.  If $L$ is strongly non-degenerate, then the quadratic
annihilator of $L$ inside $Q$ coincides with the annihilator of
$L$ in $Q$, and both are zero (Theorem \ref{Coruno} (ii)).

Another application of Theorem~\ref{qadann} leads to the proof of
the fact that, if a Lie algebra $L$ contains an essential ideal
which is strongly non-degenerate, then the algebra $L$ is itself
strongly non-degenerate (Proposition~\ref{Cordos} (ii)). This fact
was already proved by Zelmanov in \cite{Zel} by making use of the
Kostrikin radical, while our proof is based on elements.

According to the above, and in order to obtain the result
announced in the abstract (Theorem~\ref{nodegderoverk}) we need to
produce an essential ideal inside $\der (A)$ which is strongly
non-degenerate. The natural candidate for this is the ideal
$\mathrm{Inn}(A)$ of the so-called inner derivations of $A$, which
can be identified with the quotient $A/Z(A)$, known to be strongly
non-degenerate under appropriate mild hypotheses. However, this
ideal might fail to be essential, and this is somehow measured by
the ideal $I_Z$ of those derivations that map $A$ into the center
$Z(A)$. Our result then asserts that $\der(A)/I_Z$ is strongly
non-degenerate. In the particular case that the center $Z(A)$ of
$A$ does not contain associative ideals (e.g.~if $A$ is prime),
one has that $I_Z=0$, and then we do obtain that $\der(A)$ is
strongly non-degenerate.

Our arguments can be subsequently adjusted with some extra effort
to the case of a $*$-semiprime algebra $A$ and the  Lie algebra
$\sder(A)$ of those (associative) derivations of $A$ that commute
with the involution, that is, those $\delta\in \der(A)$ such that
$\delta(a^\ast)=(\delta(a))^\ast$ for every $a\in A$ (Theorem
\ref{nodegderoverkforK}).

\section{Notation and preliminaries}

Let $\Phi$ be a unital commutative ring. All algebras in this
paper, associative or not, will be $\Phi$-modules. Recall that a
\emph{Lie algebra} over $\Phi$ is a $\Phi$-module $L$, together
with a bilinear map $[\ ,\ ]\colon L\times L \to L$, denoted by
$(x, y)\mapsto [x, y]$ and called the \emph{bracket} of $x$ and
$y$ such that the following axioms are satisfied:
\begin{enumerate}[(i)]\itemsep=2mm
\item $[x, x]=0$,
\item $[x, [y, z]]+[y, [z, x]]+[z, [x, y]]=0$ (\emph{Jacobi identity}),
\end{enumerate}
\par
\noindent for every $x$, $y$, $z$ in $L$.

The standard example is obtained by considering a (non-necessary unital) associative algebra $A$, with its
same module structure and bracket given by $[x, y]=xy-yx$. Sometimes the notation $A^{-}$ is used in order to
emphasize the Lie structure of $A$.

Given an element $x$ of a Lie algebra $L$, we may define a map $ \ad x\colon L \to L$ by $\ad x(y)=[x, y]$
(which is a derivation of the Lie algebra $L$). We shall denote by $A(L)$ the associative subalgebra
(possibly without identity) of $\mathrm{End}(L)$ generated by the elements $\ad x$ for $x$ in $L$.


An element $x$ in a Lie algebra $L$ is an \emph{absolute zero divisor} if $(\ad x)^2=0$. This is equivalent
to saying that $[x, [x, L]]= 0$. The algebra $L$ is said to be \emph{strongly non-degenerate} (according to
Kostrikin) if it does not contain non-zero absolute zero divisors.

Given a Lie algebra $L$, we say that $L$ is \emph{semiprime} if we have $I^2\neq 0$ whenever $I$ is a
non-zero ideal. It is obvious from the definitions that strongly non-degenerate Lie algebras are semiprime,
but the converse does not hold (see~\cite[Remark 1.1]{Siles})

Next, $L$  is said to be {\emph prime} if $[I, J]\neq 0$ for any pair of non-zero ideals $I$, $J$ of $L$. An ideal $I$ of $L$  is said to be {\emph essential} if its intersection with any non-zero ideal is again a non-zero ideal.

For two subsets $X$, $Y$ of a (Lie) algebra $L$ we define the
\textit{annihilator of} $Y$ \textit{in} $X$ as the set
\[
\An_X(Y):=\{x \in X \ \vert \  [x, Y]=0\}\,,
\]
\noindent and the \textit{quadratic annihilator of} $Y$
\textit{in} $X$ to be the set
\[
\QAn_X(Y):=\{x \in X \ \vert \  [x, [x, Y]]=0\}\,.
\]

When $X=L$, we write $\An(Y)$ or $\An_L(Y)$ (if no confusion may
arise) and refer to it as the annihilator of $Y$. If $X=Y=L$, then
$\An{}(L)$ is called the \emph{centre} of $L$ and usually denoted
by $Z(L)$. In the case that $L=A^{-}$ for an associative algebra
$A$, then $Z(A^{-})$ agrees with the associative center $Z$ of
$A$. It is easy to check (by using the Jacobi identity) that
$\An{}(X)$ is an ideal of $L$ whenever $X$ is an ideal of $L$.
Therefore, for $A$ associative we can form the  Lie algebra
$A^-/Z$. We will be primarily interested in this type of Lie
algebras, and in Lie algebras that arise from associative algebras
with involution. If $A$ is associative and has an involution $*$,
then the set of its {\emph skew elements}
\[
K= K_A= \{ x \in A \mid  x^* = -x \}
\]
is a subalgebra of $A^-$. The center $Z(K)$ of the Lie algebra $K$ will be for brevity denoted by $Z_{K}$,
and we will be interested in the Lie algebra $K/Z_K$.

The notion of the quadratic annihilator of an (arbitrary) algebra -- defined in a similar way as we have done for Lie algebras -- plays an important role: see, for example, Smirnov's
paper~\cite{Smir}. Let us remark here that the quadratic annihilator need not be closed under sums in the
case of an associative product (for an example, see~\cite{Smir}). The same phenomenon occurs
in the Lie context, as is shown in the examples below.

\begin{ejs}\rm{(1)\ Let $F$ be any field and let $L=\mathfrak{t}
(3, F)$ be the Lie algebra of upper triangular matrices (see, e.g.~\cite{Hump}). Then
$\QAn (L)=\{ a(e_{11}+e_{22}+e_{33})+ be_{13}+ce_{23}\ \vert\ a,
b, c\in F\} \cup \{a(e_{11}+e_{22}+e_{33})+ be_{12}+ce_{13}\
\vert\ a, b, c\in F\}$, where, as usual, $e_{ij}$ denotes the
matrix in ${\mathbb{M}}_3(F)$ whose entries are all zero except
for the one in row $i$ and column $j$. Then $\QAn(L)$ is not
closed under sums.

(2)\ Now, for $L$ as in (1), consider the Lie algebra $\overline{L}: =L/Z$. Then
\[\QAn(\overline{L})=\{ \overline{ae_{13}+be_{23}}\ \vert\ a, b \in F\} \cup \{\overline{ae_{12}+be_{13}}\ \vert\ a, b\in F\}\,,
\]
where $\overline{x}$ denotes the class of an element $x$ in $L$. Again we have that the quadratic annihilator of
this algebra $\overline{L}$ is not closed under sums.}
\end{ejs}

Let $L\subseteq Q$ be Lie algebras. When $0\neq [L, q]\subseteq L$ for every non-zero $q\in Q$ we say that
$Q$ is a \emph{weak algebra of quotients of} $L$ (see~\cite{Siles}). The notion of algebra of quotients of an
algebra (associative or not necessarily associative) has a long history and is an active research area,
specially in recent years, following its development in the Lie and Jordan contexts. In the seminal
paper~\cite{Siles} the second author initiated the study of algebras of quotients of Lie algebras, by
adapting some ideas from the associative and also Jordan (\cite{mart}) contexts. She introduced the notion of
a general (abstract) algebra of quotients of a Lie algebra, and also the notion of the maximal algebra of
quotients $Q_m(L)$ of a semiprime Lie algebra $L$. Follow up results can be found in~\cite{PS, CaSa, BPSS}.


Let $B$ be a subalgebra of an associative algebra $A$. A linear
map $\delta:B\to A$ is called a {\emph derivation} if $\delta
(xy)= \delta (x)y + x \delta (y)$ for all $x, y\in B$. By a
derivation of $A$ we simply mean a derivation from $A$ into $A$.
Let $\der (A)$  denote the set of all  derivations of $A$. It is
clear that $\der (A)$ becomes a $\Phi$-module under natural
operations and it also becomes a Lie algebra by putting $[\delta,
\mu]=\delta \mu- \mu\delta$, for every $\delta$, $\mu$ in
$\der(A)$. Any element $x$ of $A$ determines a map $\ad x\colon
A\to A$ defined by $\ad x (y)=[x, y]$, which is a derivation of
$A$. For every Lie ideal $U$ of $A$, the restriction of the map
$\ad\colon A\to\der (A)$ to $U$,
$$\begin{matrix}
& U & \to & \der (A) \cr & y & \mapsto & \ad y
\end{matrix}$$
defines a Lie algebra homomorphism with kernel $\Ann_U (A)$, which allows us to identify $U/\Ann_U (A)$ with the
subalgebra $\ad (U)$ of $\der (A)$.   For any $y\in U$ and $\delta \in \der (A)$, $[\delta, \ad y]=\ad \delta
(y)$, hence $\ad (U)$ is an ideal of $\der (A)$ whenever $\delta (U) \subseteq U$ for every $\delta \in \der
(A)$. The ideal $\ad (A)$ of $\der(A)$ is usually denoted by $\inn (A)$ and its elements are called \emph{inner derivations} of $A$. Note that $A^-/Z \cong \inn (A)$.

Now let $A$ be an associative algebra with involution $\ast$. The
set
\[
\sder(A)=\{\delta\in \der(A)\mid \delta(x^*)=\delta(x)^* ~ \mbox{for all} \ x \in A \}
\]
is a Lie subalgebra of $\der(A)$. Denote by $\Ad (K)$ the set of Lie
derivations $\Ad x: A\rightarrow A$ with $x$ in $K$.

In what follows we will assume that $2$ and $3$ are invertible
elements in $\Phi$.

\section{The results}

 By an \emph{extension of Lie algebras} $L\subseteq Q$ we will mean
 that $L$ is a (Lie) subalgebra of the Lie algebra $Q$.
 Let $L\subseteq Q$ be an extension of Lie algebras and let
 $A_Q(L)$ be the associative subalgebra of $A(Q)$ generated by
 $\{\ad{x} : x \in L\}$.

For an extension $L\subseteq Q$ of Lie algebras, the condition that $\An_L(Q)=0$ ensures that the map $L\to
A(Q)$ given by $x\mapsto \ad x$ is a monomorphism of Lie algebras. Examples of extensions where $\An_L(Q)=0$
are the dense ones (see~\cite{Ca} for the definition of a dense extension and~\cite{PS} for examples).

\begin{thm}\label{qadann} Let $L\subseteq Q$ be an extension of Lie
algebras such that the map $x\mapsto \ad x$ from $L\to A(Q)$ is a
monomorphism of Lie algebras. Let $a\in \QAn_Q(L)$. Then, for
each $u \in L$ satisfying $x:=[a, u]\in L$, we have that $z:=[x,
[x, v]]$ is in $\QAn_L(L)$, for every $v\in L$.
\end{thm}

\begin{proof}

In order to ease the notation in our computations, we shall temporarily get rid
of the prefix $\ad$ and use capital letters $X$, $Y$, etc. instead
of $\ad x$, $\ad y$, etc. Because of our assumption, we shall also identify an element $x$ of $L$ with its corresponding operator $X=\ad x$ in $A(Q)$. An equation involving commutators on $L$ is then translated into the corresponding equation with capital letters and commutators in $A(Q)$.

Let $a$ be in $\QAn_Q(L)$. Then

\begin{equation} \label{uno} [a, [a, y]]=0 \text{ for every }  y\in
L\,.
\end{equation}

This implies $[A, [A, Y]]=0$ for every $Y\in \ad(L)\subseteq
A(Q)$, hence

\begin{equation} \label{dos} A^2Y+YA^2-2AYA=0 \text{ for every }
Y\in \ad(L)\,.
\end{equation}

By (\ref{uno}) we have
\begin{equation} \label{tres} A^2=0 \text{ on }  L\,,
\end{equation}

\noindent and by (\ref{dos}) and $\frac{1}{2}\in \Phi$,
\begin{equation} \label{cuatro} AYA=0 \text{ on $L$, for every }\
Y\in \ad(L).
\end{equation}

For $x=[a, u]\in L$, with $u\in L$, we have that
\begin{equation}
\label{cinco} X^2=(AU-UA)^2=AUAU- AU^2A-UA^2U+UAUA= -AU^2A
\end{equation}

 on $L$, by using (\ref{tres}) and (\ref{cuatro}).



Note that $X^3 \stackrel{(\ref{cinco})}{=} (-AU^2A)(AU-UA)=
 -AU^2A^2U+AU^2AUA=0$ by (\ref{tres}) and (\ref{cuatro}). Thus:

\begin{equation} \label{seis} X^3=0 \text{ on } L\,.
\end{equation}





We may now apply~\cite[Lemma 1.5.6]{Kos} in order to obtain:

\begin{equation} \label{amayuscula} XYX^2=X^2YX \text{ on $L$, for every }\
Y\in \ad(L)\,,\text{ and}
\end{equation}

\begin{equation} \label{cmayuscula} (YX^2)^2=X^2Y^2X^2 \text{ on $L$, for every }\
Y\in \ad(L)\,.
\end{equation}

Now,
\begin{equation}
\label{delreferee}
\begin{split}
X^2Y^2X^2 & = (YX^2)^2  \,\,\,\text{(by~(\ref{cmayuscula}))}\\
& = YX^2YX^2  \\
& = Y(XYX^2)X  \,\,\,\text{(by~(\ref{amayuscula}))}\\
& =0  \,\,\,\text{(by~(\ref{seis}))}\\
\end{split}
\end{equation}

Taking $z:=[x, [x, v]]$, with $v\in L$, we compute that, on $L$,
\begin{equation*}
\begin{split}
Z^2 & = (X^2V+VX^2-2XVX)(X^2V+VX^2-2XVX)\\
& = X^2VX^2V+VX^4V-2XVX^3V \\
&\,\,\,\,\,\, +X^2V^2X^2+VX^2VX^2-2XVXVX^2\\
& \,\,\,\,\,\, -2X^2VXVX-2VX^3VX+4XVX^2VX\\
& = X^2VX^2V+VX^2VX^2\\
& \,\,\,\,\,\, -2XV(XVX^2)-2(X^2VX)VX+4XVX^2VX \,\,\text{(by~(\ref{delreferee}) with $Y=V$ and~(\ref{seis}))}\\
& =X^2VX^2V+VX^2VX^2\,\,\text{(by~(\ref{amayuscula}) twice with $Y=V$)}\\
& =X^2VX^2V\,\,\text{(by~(\ref{delreferee}) with $Y=V$)}\\
& =(X^2VX)XV=(XVX^2)XV\,\,\text{(by~(\ref{amayuscula}) with $Y=V$)}\\
& =0\,,\,\,\text{(by~(\ref{seis}))}\\
\end{split}
\end{equation*}
which completes the proof of the theorem.
\end{proof}

\begin{thm}\label{Coruno} Let $L\subseteq Q$  be an extension of Lie
algebras with $Q$ a \wqa of $L$ and $L$ strongly non-degenerate.
Then:
\begin{enumerate}[{\rm (i)}]
\item $Q$ is strongly non-degenerate {\rm (\cite[Proposition
2.7(iii)]{Siles})}.
\item $\An_Q(L)=\QAn_Q(L)=0$.
\end{enumerate}

\end{thm}

\begin{proof} (i). Suppose that there exists a non-zero element $a\in
\QAn_Q(Q)$, and choose $u\in L$ satisfying $0\neq x:=[a, u]\in L$. Since $L$ is strongly non-degenerate,
$z:=[x, [x, v]]\neq 0$ for some $v\in L$. But, by  Theorem~\ref{qadann}, $z$ must be zero, a contradiction.

(ii). $\An_Q(L)=0$ because $Q$ is a \wqa of $L$. For $a\in \QAn_Q(L)$, $z:=[x,[x, v]]\in \QAn_L(L)$ whenever
$x:=[a, u]\in L$, with $u\in L$ (Theorem~\ref{qadann}). Now $L$ being strongly non-degenerete implies
$\QAn_L(L)=0$, whence $x$ is zero. Since $Q$ is a \wqa of $L$, we obtain that $a=0$.
\end{proof}

Condition (ii) in the result below was proved by Zelmanov
in~\cite[Corollary 2 in pg. 543]{Zel} for strongly non-degenerate
Lie algebras by using the Kostrikin radical. Our proof here is
based on elements.

\begin{prop}\label{Cordos} Let $I$ be a strongly non-degenerate ideal of a
Lie algebra $L$. Then:
\begin{enumerate}[{\rm (i)}]
\item $\An_L(I)=\QAn_L(I)$.
\item If $\An_L(I)=0$, then the algebra $L$ is strongly non-degenerate.
\end{enumerate}

\end{prop}

\begin{proof} (i). Clearly, $\An_L(I)\subseteq \QAn_L(I)$.
Conversely, consider $a\in \QAn_L(I)$. The strongly non-degeneracy assumption on $I$ implies that the map
$I\to A(L)$ given by $y\mapsto \ad y$ is a monomorphism of Lie algebras. Since $I$ is a strongly
non-degenerate ideal of $L$, Theorem~\ref{qadann} implies that for every $u\in I$, the element $x=[a, u]$ is
in $\QAn_I(I)=0$, hence $a\in \An_L(I)$.

(ii). In this case we have that $L$ is a \wqa of $I$. Apply
Theorem~\ref{Coruno} to obtain that $L$ must be strongly
non-degenerate too.
\end{proof}

If $A$ is an associative algebra, and since every derivation maps $Z$ to $Z$, the set
\[
I_Z=\{\delta\in\der (A)\mid \delta(A)\subseteq Z\}
\]
is easily seen to be a Lie ideal of $\der (A)$ that contains the center of $\der (A)$. Indeed, for every
$\delta\in Z(\der(A))$ and each $a\in A$, we have that $0=[\delta, \ad a]= \ad (\delta(a))$, hence $\delta(a)\in Z$.
Moreover, under certain conditions, $\mathrm{Inn}(A)$ can be seen as an essential ideal of $\der(A)/I_Z$.

\begin{lem}\label{I_Z} Let $A$ be a semiprime non-commutative associative algebra. Then:
\begin{enumerate}[{\rm (i)}]
\item $\mathrm{Inn}(A)$ is (isomorphic to) an essential ideal of $\der(A)/I_Z$, where $I_Z$ is defined as
before.
\item If $Z$ does not contain non-zero
associative ideals (in particular, if $A$ is prime), then $I_Z=0$.
\end{enumerate}

\end{lem}

\begin{proof} (i). The map
\[
\begin{matrix} \mathrm{Inn} (A) & \to & \der(A)/I_Z \cr \ad a & \mapsto &
\overline{\ad a}

\end{matrix}
\]
is a monomorphism of Lie algebras.
This follows from the following fact:
\[
\label{doble}
 \tag{\ensuremath{\dagger}} [[a, A], A]=0\,,\text{ with } a\in A,  \text{ implies }\  a\in Z\,.
\]
Indeed, $[[a, A], A]=0$ implies, by~\cite[Sublemma in pg. 5]{HTopics}, $[a, A]=0$, that is, $\ad a=0$.

This allows us to identify $\mathrm{Inn}(A)$ with its image inside $\der(A)/I_Z$. The formula $[\delta,\ad a]=\ad \delta(a)$, where $a\in A$ and $\delta\in\der(A)$ ensures that $\mathrm{Inn}(A)$ is indeed an ideal of $\der(A)/I_Z$.

Now let $J/I_Z$ be a non-zero ideal of $\der(A)$ and consider
$\delta\in J \setminus I_Z$, that is, $[\delta(A),A]\neq 0$. A
second usage of~(\ref{doble}) allows us to conclude that
$[[\delta(A),A], A]\neq 0$. Take $a$ in $A$ such that
$[\delta(a),A]\not\subseteq Z$. Then $\overline 0 \neq
\overline{[\delta, \ad a]}=\overline{\ad \delta(a)}$, and thus
$\mathrm{Inn}(A)$ is essential in $\der(A)/I_Z$.

(ii). Take $\delta\in I_Z$ and $d\in \der(A)$. Put $\mu=[\delta,
d]$. For every pair of elements $a$, $b\in A$ we have $\mu([a,
b])=[\mu(a), b]+[a, \mu(b)]$. Note that $[\mu(a), b]=[\delta d(a),
b]-[d\delta (a), b]=-[d\delta (a), b]=0$, because $d(Z)\subseteq
Z$, and analogously $[a, \mu(b)]=0$. It follows from this that
$\mu([A, A])=0$.

Now let $I$ be a non-central Lie ideal of $A$, and take $y$ in $I\setminus Z$. Then $[y, A]\neq 0$ and by (\ref{doble}), we get $0\neq [[y, A],A]\subseteq I\cap [A, A]$. Thus $[A, A]$ intersects non-trivially every non-central Lie ideal of $A$.

We claim that the subalgebra $\la [A, A]\ra$ generated by $[A,A]$ contains an essential associative
ideal of $A$. Herstein's~\cite[Theorem 3]{hers70} implies that $\la [A, A]\ra$ contains a non-zero associative ideal. By Zorn's Lemma, it is possible to find $M$ maximal among all the associative ideals contained in $\la [A, A]\ra$. If $\Ann(M)$ were non-zero, we get by what we have just proved that $\Ann(M)\cap [A, A]\neq 0$.
Again by~\cite[Theorem 3]{hers70}, and since $Z$ does not contain non-zero associative ideals, we have that $\la \Ann(M)\cap [A, A]\ra$ contains a non-zero associative ideal $J$. Notice that $J$ is not contained in $M$, since otherwise $J\subseteq M\cap\Ann (M)$, which is zero because $M$ is non-zero and $A$ is semiprime. Then $M\oplus J\subseteq M\oplus \Ann(M)\cap \la [A, A]\ra$ and $M\subsetneq M\oplus J$, which contradicts the maximality of $M$.

By the first part of the proof we have $\mu([A, A])=0$, and so
$\mu(I)=0$, where $I$ is an essential ideal of $A$ contained in
$\la [A, A]\ra$. This implies that $\mu=0$. For, if $\mu(a)\neq 0$
for some $a\in A$, then the essentiality of $I$ implies that there
exists $y\in I$ such that $y\mu(a)\neq 0$. But
$y\mu(a)=\mu(ya)-\mu(y)a=0$, a contradiction.
\end{proof}

\begin{thm}\label{nodegderoverk}

Let $A$ be a semiprime non-commutative associative algebra. Then:
\begin{enumerate}[{\rm (i)}]
\item $\der(A)/I_Z$ is a strongly non-degenerate Lie algebra.
\item If $Z$ does not contain non-zero associative
ideals (in particular, if $A$ is prime), then $\der(A)$ is a strongly non-degenerate Lie algebra.
\end{enumerate}

\end{thm}
\begin{proof} (i). Use~\cite[Lemma 5.2]{DFGG}, Proposition~\ref{Cordos} (ii) and Lemma~\ref{I_Z} (i).

(ii) follows from (i) and condition (ii) in Lemma~\ref{I_Z}.
\end{proof}

We now consider the case where our associative algebra $A$ has an
involution $*$. Under some additional mild assumptions on $A$ in
order to rule out algebras of low degrees, we obtain similar
results on the non-degeneracy of $\sder(A)$. Rather than proving
them in full, we just indicate which changes are needed to adjust
Lemma~\ref{I_Z} and Theorem~\ref{nodegderoverk} to the current
setting.

Recall that, if $A$ is a semiprime associative algebra, the \emph{extended centroid} $\mathcal{C}=\mathcal{C}(A)$ of $A$ is defined as the center of the two-sided right ring of quotients. It also coincides with the center of $Q_s(A)$, the symmetric ring of quotients. For every $x$ in an algebra $A$ we define $\deg(x)$ as the degree
of algebraicity of $x$ over the extended centroid $\mathcal{C}$, provided that $x$ is algebraic. If $x$ is not algebraic,
then we define $\deg(x)=\infty$. Put $\deg(A)=\sup\{ \deg(x)\mid
x\in A\}$. It is well-known that $\deg(A) < \infty$ if and only if
$A$ is a PI algebra. Furthermore, it is known that $\deg(A) = n <
\infty$ if and only if $A$ satisfies the standard polynomial
identity of degree $2n$, but does not satisfy any polynomial
identity of degree $< 2n$, and this is further equivalent to the
condition that $A$ can be embedded into the matrix algebra
$\mathbb{M}_n(F)$ for some field $F$ (one can take, say, $F$ to be
the algebraic closure of $\mathcal{C}$), but cannot be embedded
into $\mathbb{M}_{n-1}(R)$ for any commutative algebra $R$.

\begin{enumerate}[{\rm (a)}]
\item
The analogue of (\ref{doble}) in Lemma~\ref{I_Z} is as follows:

If $A$ is a semiprime algebra with involution, then if $a$ belongs to the skew elements $K$ in $A$
and $[[a, K], K]=0$ we get $[a, K]=0$. To prove this, suppose $a
\in K$ such that $[a, K]\neq 0$. This means that $\overline a\neq
0$, where $\overline a$ denotes the class of $a$ in $K/Z_K$. But
then $[\overline{a}, K/Z_K]\neq 0$, since $K/Z_K$ is semiprime
by~\cite[Theorem 3]{Jordan80}, that is, $[[a, K], K]\neq 0$.
\item
The use of~\cite[Theorem 3]{hers70} in the proof of condition (ii)
in Lemma ~\ref{I_Z}  must be changed to~\cite[Lemmas 2 and
3]{lans76}. In order to apply these results, certain restrictions
on the degree of the algebra are needed. Concretely,
$\mathrm{deg}(A/I) > 2$ for every $\ast$-prime ideal $I$ of $A$.

\end{enumerate}

Recall that an ideal $I$ in an algebra $A$ with involution $*$ is
a \emph{$\ast$-ideal} if $I$ is invariant under the involution,
that is, $I^\ast=I$. The algebra $A$ is said to be
$\ast$-\emph{prime} if the product of two non-zero $*$-ideals is
again non-zero. A $*$-ideal $I$ is said to be $\ast$-\emph{prime}

if $A/I$ is a $*$-prime algebra. The definition of a
$*$-\emph{semiprime} algebra is analogous.

Define the following Lie ideal of $\sder(A)$:
\[
I_{K,Z}=\{\delta\in\sder (A)\mid \delta(K)\subseteq Z\}\,.
\]

In the current context, our Lemma~\ref{I_Z} takes then the following form.

\begin{lem}\label{I_ZparaK} Let $A$ be a semiprime
non-commutative associative algebra with involution $\ast$. Then:
\begin{enumerate}[{\rm (i)}]
\item $\mathrm{Inn}(K)$ is (isomorphic to) an essential ideal of
$\sder(A)/I_{K,Z}$, where $I_{K,Z}$ is defined as before.
\item If $Z(K)$ does not contain non-zero associative
 $\ast$-ideals (in particular, if $A$ is $\ast$-prime), then $I_{K,Z}=0$.
\end{enumerate}

\end{lem}

The analogue of~\cite[Lemma 5.2]{DFGG} (used in the proof of condition (i) in Theorem~\ref{nodegderoverk}) is the proposition below,
which again requires conditions on the degree of the algebra. In
particular, it generalizes ~\cite[Theorem 2.13]{Benkart}. Recall
that an involution $\ast$ in an associative algebra $A$ is said to
be of the \emph{first kind} if it is the identity on the centroid
of $A$. Otherwise it is called of the \emph{second kind}.

\begin{prop}\label{SND} Let $A$ be a $*$-semiprime algebra. Assume that the involution is either of the second kind, or else it is of the first kind and $\mathrm{deg}(A/I)> 2$ for every $\ast$-prime ideal $I$ of $A$. Then $[k, [k,
K]]\subseteq Z(A)$, with $k\in K$, implies $k\in Z(A)$. In
particular, $K/{(K\cap Z)}$ is a strongly non-degenerate Lie
algebra.
\end{prop}

\begin{proof}

Let $I$ be a $\ast$-ideal of $A$. It is clear that $A/I$ also
becomes a $*$-algebra with the natural involution.

On the other hand, if $\overline{x}$ denotes the class of an
element $x$ in $A/I$ and $\overline{K}=\{\overline{k}\mid k\in
K\}$, we have that $\overline{K}=K_{A/I}$. The containment
$\overline{K}\subseteq K_{A/I}$ is clear, and for the converse,
take $\overline{a}$ in $K_{A/I}$ and let $y\in I$ be such that
$a^\ast+a=y$. Then $(a^\ast -
{\frac{1}{2}}y)^\ast=y-a-{\frac{1}{2}}y = -a+ {\frac{1}{2}}y$,
that is, $a-{\frac{1}{2}}y\in K$, and
$\overline{a}=\overline{a-{\frac{1}{2}}y}$.

Now, consider $k\in K$ satisfying $[k, [k, K]]\subseteq Z$. In
particular, $({\ad k})^3 (t)=0$ for every $t\in K$. Arguing as in
the proof of~\cite[Lemma 5.2]{DFGG}, we obtain:
\[
\label{triple} \tag{\ensuremath{\ddagger}}  (\ad k)^2(t)=0\,.
\]

Let $\{I_\alpha\}_{\alpha\in \Lambda}$ be the collection of all
$\ast$-prime ideals of $A$. Since $A$ is $\ast$-semiprime,
$\bigcap_{\alpha\in \Lambda}I_\alpha=0$.  Suppose $[\overline{k},
A/I_\beta]\neq 0$ for some $\beta\in \Lambda$. Since $A/I_{\beta}$
is a $*$-prime algebra we may apply~\cite[Lemma 5.4]{BPSS} in
order to conclude that $Z(\overline{K})=Z(A/I_\beta)\cap
\overline{K}$, and hence $[\overline{k}, \overline{K}]\neq 0$.
Use~\cite[Theorem 5.3]{BPSS} if $\ast\colon A/I_\beta\to
A/I_\beta$ is of the first kind, or~\cite[Theorem 2.13]{Benkart}

if the involution is of the second kind, to conclude that
$[\overline{k}, [\overline{k}, \overline{K}]]\neq 0$, in
contradiction to (\ref{triple}). In consequence, $[\overline{k},
A/I_\alpha]= 0$ for every $\alpha \in \Lambda$, that is, $[k,
A]\subseteq \bigcap_{\alpha\in \Lambda}I_\alpha =0$.
\end{proof}

Finally, the involutive version of Theorem~\ref{nodegderoverk} is the following:

\begin{thm}\label{nodegderoverkforK}

Let $A$ be a $*$-semiprime non-commutative associative algebra
with $\mathrm{deg}(A/I) > 2$ for every $\ast$-prime ideal $I$ of
$A$. Then:
\begin{enumerate}[{\rm (i)}]
\item $\sder(A)/I_{K,Z}$ is a strongly non-degenerate Lie algebra.
\item If $Z(K)$ does not contain non-zero associative
$\ast$-ideals (in particular, if $A$ is $\ast$-prime), then
$\sder(A)$ is a strongly non-degenerate Lie algebra.
\end{enumerate}

\end{thm}

\section*{acknowledgements}

Part of this work was carried out during visits of the
first author to the Universidad de M\'alaga, and of the second one to the Centre de Recerca Matem\`{a}tica. We wish to thank the host centers for their warm hospitality. We also thank the referee for a careful reading of the manuscript and for simplifying our initial proof of Theorem~\ref{qadann}.



\begin{thebibliography}{99}


\bibitem{Benkart}
\textsc{G.~Benkart}, The Lie in

ner ideal structure of associative rings, \emph{J. Algebra} \textbf{43} (1976),  561--584.

\bibitem{BPSS}


\textsc{M. Bre\v sar, F. Perera, J. S\'anchez Ortega, M. Siles Molina}, Computing the
maximal algebra of quotients of a Lie algebra. (Preprint).

\bibitem{Ca}


\textsc{M.~Cabrera}, Ideals which memorize the extended centroid,
\emph{J. Algebra Appl.} \textbf{1} (2002), 281--288.


\bibitem{CaSa}


\textsc{M.~Cabrera, J.~S\'anchez~Ortega}, Lie quotients for skew Lie algebras. (Preprint).

\bibitem{DFGG} \textsc{C. Draper Fontanals, A. Fern\'andez L\'opez, E.
Garc\'\i a, M. G\'omez Lozano}, The socle of a nondegenerate Lie
algebra, \emph{J. Algebra} \textbf{280} {(2004)}, 635--654.

\bibitem{HTopics} \textsc{I.~N.~Herstein}, \emph{Topics in Ring Theory}.
The University of Chicago Press (1969).

\bibitem{hers70} \textsc{I.~N.~Herstein}, On the Lie structure of an
associative ring, \emph{J. Algebra} \textbf{14} (1970), 561--571.

\bibitem{Hump} \textsc{ J. E. Humphreys}, \emph{Introduction to Lie Algebras
and Representation Theory}. Springer-Verlag New York Inc., 1972.

\bibitem{Jordan80} \textsc{D.~A.~Jordan}, The Lie ring of symmetric derivations of a ring with involution,
 \emph{J. Austral Math. Soc.} \textbf{29} (1980),
153--161.

\bibitem{JJ} \textsc{C.~R.~Jordan, D.~A.~Jordan}, Lie rings of derivations
of associative rings, \emph{J. London Math. Soc.} \textbf{17} (1978), 33--41.

\bibitem{Kos} \textsc{A. I. Kostrikin}, \emph{Around Burnside}.
Springer-Verlag Berlin Heidelberg, 1990.

\bibitem{lans76}\textsc{C.~Lanski}, Lie structure in semi-prime rings
with involution, \emph{Comm. Alg.} \textbf{4} (1976), 731--746.

\bibitem{mart} \textsc{C. Mart\'\i nez}, The ring of fractions of a
Jordan algebra, \emph{J. Algebra} \textbf{237} (2001), 798--812.

\bibitem{PS} \textsc{F. Perera, M. Siles Molina}, Associative and Lie algebras of quotients, \emph{Publ. Mat.}
(To appear.)

\bibitem{Siles} \textsc{M. Siles Molina}, Algebras of quotients of Lie
algebras, \emph{J. Pure Appl. Algebra} \textbf{188} {(2004)},
175--188.

\bibitem{Smir} \textsc{O. Smirnov}, Finite $\mathbb{Z}$-gradings of Lie
algebras and symplectic involutions, \emph{J. Algebra} \textbf{218} {(1999)},
246--275.

\bibitem{Zel} \textsc{E. Zelmanov}, Lie algebras with an algebraic adjoint
representation. \emph{Math. USSR Sb.} {\bf 49}(2) {(1984)}, 537--552.

\end{thebibliography}
\end{document}